\documentclass[12pt]{amsart}

\usepackage{amsthm}
\usepackage{amscd}
\usepackage{amsmath}
\usepackage{amsfonts}
\usepackage{amssymb}
\usepackage{fullpage}
\usepackage{graphicx}

\def\hide#1{}

\newenvironment{pf}{
\medskip
  \small
  \Pf
}{
  \hfill $\bigcirc$
  \normalsize
}

\def\Pf{\medskip\emph{Pf}. }

\def\bb{\mathbb}

\makeatletter

\def\underrightharpoondown#1{\mathop{\vtop{\m@th
 \ialign{##\crcr
 $\hfil\displaystyle{#1}\hfil$\crcr
 \noalign{\kern0.5pt\nointerlineskip}\rightharpoondownfill\crcr  
 \noalign{\kern0.5pt}}}}}

\def\rightharpoondownfill{$\m@th\smash-\mkern-7mu
 \cleaders\hbox{$\mkern-2mu\smash-\mkern-2mu$}\hfill
 \mkern-7mu\mathord\rightharpoondown$}

\makeatother

\def\hat{\widehat}

\def\hat{\widehat}

\def\theorem{\medskip\noindent\textbf{Theorem. }}

\def\corollary{\medskip\noindent\textbf{Corollary. }}

\newcommand{\vol}[1]{{\bf #1}}
 
\def\lb{\langle}
\def\rb{\rangle}

\def\e{\textbf{e}}

\let\hat\widehat

\def\and{\quad\text{and}\quad}

\begin{document}  

\title[Notes on the equivalence of first-order rigidity]{
	Some notes on the equivalence of first-order 
	rigidity in various geometries
	}
\author[F. Saliola]{
	Franco V. Saliola$^\dag$
	}
\thanks{
	$^\dag$Department of Mathematics, 310 Malott Hall,
	Cornell University, Ithaca, NY 14853-4201 USA
	(saliola@math.cornell.edu)
	Work supported, in part, by an NSERC (Canada) research award
	held at York University and an NSERC (Canada) PGS A award.
	}
\author[W. Whiteley]{
	Walter J. Whiteley$^\ddag$
	}
\thanks{
	$^\ddag$Department of Mathematics and Statistics,
	York University, 4700 Keele Street,
	North York, Ontario, M3J 1P3 Canada
	(whiteley@mathstat.yorku.ca).
	Work supported by a grant from NSERC (Canada).
	}
 
\begin{abstract} 
These pages serve two purposes. First, they are notes to accompany the
talk \textit{\textbf{Hyperbolic and projective geometry in constraint
programming for CAD}} by Walter Whiteley at the \emph{J\'anos Bolyai
Conference on Hyperbolic Geometry}, 8--12 July $2002$, in Budapest,
Hungary. Second, they sketch results that will be
included in a forthcoming paper that will present the equivalence of
the first-order rigidity theories of bar-and-joint frameworks in
various geometries, including Euclidean, hyperbolic and spherical
geometry. The bulk of the theory is outlined here, with remarks and
comments alluding to other results that will make the final version of
the paper. 
\end{abstract}

\maketitle


\section{Introduction}

In this paper, we explore the connections among 
the theories of first-order rigidity of bar and joint
frameworks (and associated structures) in various metric
geometries extracted from the underlying projective space of
dimension $n$, or
$\bb R^{n+1}.$  The standard examples include Euclidean space,
elliptical (or spherical) space, hyperbolic space, and a metric
on the exterior of hyperbolic space.

In his book, Pogorelov explored more general issues of
uniqueness, and local uniqueness of realizations in these
standard spaces, with some first-order correspondences as
corollaries \cite{pogo}.  
We will take the opposite tack --
beginning directly with the first-order theory, in this paper. 
We believe this presents a more transparent and accessible
starting point for the correspondences.  In a second paper, we
will use the additional technique of `averaging' in combination
with the first-order results to transfer results about pairs of
objects with identical distance constraints in one space to
corresponding pairs in a second space 
%
%
%
%

Like Pogorelov (and perhaps for related reasons) we will begin
with the correspondence between the theory in
elliptical or spherical space and the theory in
Euclidean space (\S \ref{equivalence S -> E}).  
This correspondence of configurations
is direct -- using gnomic projection (or central projection) from
the upper half sphere to the corresponding Euclidean space. 
This correspondence between spherical frameworks and their
central projections into the plane is also embedded in previous
studies of frameworks in dimension $d$ and their one point cones
into dimension $d+1$ \cite{Wh2}.  

With a firm grounding for the first-order rigidity in
spherical space, it is simpler to work from the spherical
$n$-space to the other metrics extracted from the underlying
$\bb R^{n+1}$ (\S\ref{equivalence in other geometries}).
The correspondence works for any metric of
the form $\lb p,q \rb = \sum_{i=1}^{n+1} a_{i} p_i q_i$, $a_{i} \neq 0$,
in addition to the special case of Euclidean space (with
$a_{n+1}=0$).   It has a particularly simple form, for selected
normalizations of the rays as points in the space, such as
$\lb p,p \rb = \pm 1$, which is the form we present.  

Having examined the theory of first-order motions, we pause to
present the motions as the solutions to a matrix equation 
$R_X (G,p) x = 0$ for the metric space $X$ (\S\ref{rigidity matrix}).
In this
setting, we have the equivalent theory of static rigidity working
with the row space and row dependences (the self-stresses) of
these matrices, instead of the column dependencies (the
motions).   The correspondence is immediate, but it takes a
particular nice form for the `projective' models in Euclidean
space of the standard metrics.  In this setting, the rigidity
correspondence is a simple matrix multiplication:
	$$
	R_X (G,p) [T_{XY}] =  R_Y (G,p)
	$$
for the same underlying configuration $p$, where $[T_{XY}]$ is
a block diagonal matrix with a block entry for each vertex,
based on how the sense of `perpendicular' is twisted at that
location from one metric to the other.
As a consequence of this simple correspondence of matrices,
we see that row dependencies (the static self-stresses) are
completely unchanged by the switch in metric.  As a biproduct
of this static correspondence, there is a correspondence for
the first-order rigidity of the structures with inequalities,
the tensegrity frameworks, which are well understood as a
combination of first-order theory and self-stresses of the
appropriate signs for the edges with pre-assigned
inequality constraints.  

As this shared underlying statics hints, there is a
shared underlying projective theory of statics
(and associated first-order kinematics) 
\cite{CrapoWhiteley}.

We will not present that theory here but we note
the projective invariance, in all the metrics, of the
first-order and static theories 
%
%
%
%
	(\S7).  
There are various
extensions that follow from this underlying projective
theory, such as inclusion of `vertices at infinity' in
Euclidean space \cite{CrapoWhiteley},
and the
possibility that polarity has a role to play (see below).  

As an application of these correspondences, we consider a
classical theory of rigidity for polyhedra -- the theorems of
Cauchy, Alexandrov, and the associated theory of Andreev. 
This theory provides theorems about the first-order
rigidity of convex polyhedra and convex polytopes with
either rigid faces, or $2$-faces triangulated with bars and
joints in dimensions
$d\geq 3$, in Euclidean space.  Since the basic concepts of
convexity transfer among the metrics (if we remove the equator
on the sphere, or the corresponding line at infinity in
Euclidean space), this first-order and static theory immediately
transfers to identical theorems in the other metric spaces
(\S\ref{andreev}).  There are some first-order extensions of Cauchy's
Theorem to versions of local convexity, which will automatically
extend to the various metrics and on through to hyperplanes and
angles, giving additional generalizations.
%
%
%
%
%
%
%
%
%
%
 Moreover,
this theory for hyperplanes and angles will be projectively
invariant, if we are careful with the transfer of concepts such
as `convexity' through the projective transformations. 

In hyperbolic space, there is a correspondence between rigidity
of `bar-and-joint frameworks' with vertices and distance constraints
in the exterior hyperbolic space (or ideal points) and planes
and angle constraints in the interior hyperbolic space.  
We present this correspondence directly, although it can be
viewed as a polarity about the absolute. 
With this correspondence, the first-order Cauchy theory in
exterior hyperbolic space gives a first-order theory for planes
and angles in hyperbolic space.  This result turns out to be a 
generalization of the
first-order version of Andreev's Theorem. In this setting, the
constraint that angles be less than $\pi/2$ disappears and the
angles have the full range of angles in a convex polyhedron ($<
\pi$).     

Moreover, as this hints, there is a correspondence, via
spherical polarity, which connects the first-order Cauchy
Theorem in the spherical or elliptic space with an
Andreev style first-order theorem for planes and angles of a
simple convex polytope in elliptical geometry (\S\textbf{none}).  
The effect of polarity in Euclidean space is drastically different.
It has an interesting, and distinctive interpretations in
dimensions $d=2$ and $d=3$ \cite{Wh5, Wh6}.  

%
%
%
%

The general problem of characterizing which graphs have some
(almost all) realizations in $d$-space as first-order rigid
frameworks is hard for dimensions $d\geq 3$.  With these
correspondences, we realize that this problem is identical in
all the metric spaces and we will not get additional leverage
by comparing first-order behaviour under the various metrics. 

On the other hand, in general geometric constraint programming
in fields such as CAD, there is an interest in more general
systems of geometric objects and general constraints.  For
example, circles of variable radii with angles of intersection
as constraints are in interest in CAD.  As people familiar
with hyperbolic geometry may realize, these are equivalent,
both a first-order and at all orders, to planes and angles in
hyperbolic 3-space.  The correspondence presented here provides
the final step in the correspondence between circles and angles
in the plane and points and distances in Euclidean $3$-space
\cite{SaliolaWhiteley}.

The basic first-order correspondence among metrics should
extend to differentiable surfaces from these discrete
structures.  The major difference here is that static rigidity
and first-order rigidity are distinct concepts in the this
world which corresponds to infinite matrices.  Still the
correspondence should apply to both theories, and all the
metrics.  


\section{First-Order Rigidity in $\bb E^n$}

\subsection{Euclidean $n$-space} Let $\bb E^n$ denote the set of vectors
in $\bb R^{n+1}$ with $x_{n+1} = 1$,
	$$\bb E^n = \{ x \in \bb R^{n+1} \mid \e \cdot x = 1 \},$$
where $\e = (0,0, \ldots, 1) \in \bb R^{n+1}$. 
An $m$-plane of $\bb E^n$ is the intersection of $\bb E^n$ with an 
$(m+1)$-subspace of $\bb R^{n+1}$. The distance between $x, y \in \bb E^n$ is
$d_{\bb E}(x,y) = |x - y| = \sqrt{\sum_i^n (x_i - y_i)^2}$.

\subsection{Frameworks and rigidity in $\bb E^n$}
A graph $G = (V, E)$ consists of a finite vertex set 
$V = \{1, 2, \ldots, v\}$ and an edge set $E$, 
where $E$ is a collection of unordered pairs of 
vertices. A \emph{bar-and-joint framework} $G(p)$ in $\bb E^n$ is a
graph $G$ together with a map $p: V \to \bb E^n$. 
Let $p_i$ denote $p(i)$. 

A \emph{motion} of the framework $G(p)$ is a continuous family of
functions $p(t): V \rightarrow \bb E^n$ with $p(0) = p$ such that for $\{i,
j\} \in E$, $d_{\bb E}(p_i(t), p_j(t)) = c_{ij}$, where $c_{ij}$ is a
constant, for all $t$. A framework is \emph{rigid} if all motions are
\emph{trivial}: for each $t$, there is a rigid motion $A_t$ of $\bb E^n$,
such that $A_t(p_i) = p_i(t)$, for all $i \in V$.

\subsection{Motivation for first-order rigidity.} \label{motivationE}
Suppose $p(t)$ is a motion of the framework $G(p)$ in $\bb E^n$ 
differentiable at $t = 0$. 
Since $d_{\bb E}(p_i(t), p_j(t)) = c_{ij}$ for each $\{i,j\} \in E$, 
the derivative of $p(t)$ must satisfy
	$$(p_i - p_j) \cdot (p'_i(0) - p'_j(0)) = 0,$$
where $x \cdot y$ denotes the Euclidean inner product of the vectors
$x$ and $y$. Since the framework lies in $\bb E^n$ during the motion 
($p_k(t) \in \bb E^n$ for all $k \in V$), $p_k(t)$ satisfies
$\e \cdot p_k(t) = 0$ for all $k \in V$. Hence its derivative satisfies,
	$$\e \cdot p'_i(0) = 0$$
for each $i \in V$. This motivates the following definition.

\subsection{First-order rigidity in $\bb E^n$.} \label{f.o.m in E}
A \emph{first-order motion}
of the framework $G(p)$ in $\bb E^n$ is a map $u: V \to \bb R^{n+1}$
satisfying, for each $\{i,j\} \in E$ and $k \in V$,
	\begin{equation} \label{equations for Euclidean f.o.m.}
		(p_i - p_j) \cdot (u_i - u_j) = 0 \and
		\e \cdot  u_k = 0,
	\end{equation}
where $u_i$ denotes $u(i)$. 
\vspace{-2em}
\begin{figure}[htb]
  \begin{center}
    \includegraphics[scale=0.90]{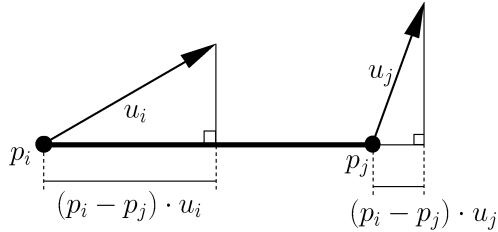}
    \caption{$u$ is a first-order motion if $(p_i - p_j) \cdot u_i =
    (p_i - p_j) \cdot u_j$ for all edges $\{i,j\}$. That is, the
    projection of $u_i$ onto $p_i - p_j$ must equal the projection of
    $u_j$ onto $p_i - p_j$.}
    \label{figure: first-order motion}
  \end{center}
\end{figure}

A \emph{trivial first-order motion} 
of $\bb E^n$ is a map $u: \bb E^n \to \bb R^{n+1}$ satisfying
	\begin{equation*}
		(x - y) \cdot (u(x) - u(y)) = 0 \and
		\e \cdot  u(z) = 0,
	\end{equation*}
for all $x$, $y$ and $z$ in $\bb E^n$.
	%
	%
$G(p)$ is \emph{first-order
rigid} in $\bb E^n$ if all the first-order motions of the framework 
$G(p)$ are restrictions of trivial first-order motions of $\bb E^n$.

\subsection{Remark} Any rigid motion of $\bb E^n$ yeilds a trivial
first-order motion of a given framework: the isometry restricts to a 
motion of the framework whose derivative satisfies the equations in 
(\ref{equations for Euclidean f.o.m.}). 

\subsection{Remark} First-order rigidity is a good indicator of
rigidity: first-order rigidity implies rigidity, but not conversely.


\section{First-Order Rigidity in $\bb S_+^n$}

\subsection{Spherical $n$-Space.} Let $\bb S_+^n$ denote the upper hemisphere
of the unit sphere in $\bb R^{n+1}$,
\[
\bb S_+^n = \{ x \in \bb R^{n+1} \mid x \cdot x = 1, \e \cdot x > 0 \},
\]
An $m$-plane of $\bb S_+^n$ is the intersection of $\bb S_+^n$ with an 
$(m+1)$-subspace of $\bb R^{n+1}$. The distance between two points
$x, y \in \bb S_+^n$ is given by the angle subtended by the vectors 
$x$ and $y$, $d_{\bb S_+}(x,y) = \text{arccos}(x \cdot y)$.

\subsection{Frameworks and rigidity in $\bb S_+^n$}
A \emph{bar-and-joint framework} $G(p)$ in $\bb S_+^n$ is a
graph $G$ together with a map $p: V \to \bb S_+^n$. 
A \emph{motion} of the framework $G(p)$ in $\bb S_+^n$ is a continuous family of
functions $p(t): V \rightarrow \bb S_+^n$ with $p(0) = p$ such that for $\{i,
j\} \in E$, $d_{\bb S_+}(p_i(t), p_j(t)) = c_{ij}$, where $c_{ij}$ is a
constant, for all $t$. A framework is \emph{rigid} if all motions are
\emph{trivial}: for each $t$, there is a rigid motion $A_t$ of $\bb S_+^n$,
such that $A_t(p_i) = p_i(t)$, for all $i \in V$.

\subsection{Motivation for first-order rigidity in $\bb S_+^n$.}
To extend the definitions of first-order motion and first-order rigidity 
to frameworks in $\bb S_+^n$, mimic the motivation presented in section 
\ref{motivationE}. If $p(t)$ is a motion of a framework $G(p)$ in 
$\bb S_+^n$, then for all $t$ and $\{i,j\} \in E$, 
	$$d_{\bb S_+}(p_i(t) \cdot p_j(t)) = c_{ij},$$
where $c_{ij}$ is constant for all $\{i,j\} \in E$,
and for all $t$ and $k \in V$,
	$$p_k(t) \cdot p_k(t) = 1.$$
Equivalently, for all $t$, $\{i,j\} \in E$ and $k \in V$,
	$$p_i(t) \cdot p_j(t) = \cos c_{ij},$$ 
	$$p_k(t) \cdot p_k(t) = 1.$$
If the motion $p(t)$ is differentiable at $t=0$, then $p(t)$ must satisfy,
	$$p_i \cdot p'_j(0) + p'_i(0) \cdot p_j = 0,$$
	$$p_k \cdot p'_k(0) = 0.$$
This leads to the following definition.

\subsection{First-Order Rigidity in $\bb S_+^n$.} A \emph{first-order motion}
of the framework $G(p)$ in $\bb S_+^n$ is a map $u: V \to \bb R^{n+1}$
satisfying, for each $\{i,j\} \in E$ and for each $k \in V$,
	\begin{equation} \label{f.o.m in S}
		p_i \cdot u_j + p_j \cdot u_i = 0
		\and
		p_k \cdot u_k = 0.
	\end{equation}	

A \emph{trivial first-order motion} of $\bb S_+^n$ is 
a map $u: \bb S_+^n \to \bb R^{n+1}$ satisfying
	\begin{equation*}
		x \cdot u(y) + y \cdot u(x) = 0 \and
		z \cdot  u(z) = 0,
	\end{equation*}
for all $x$, $y$ and $z$ in $\bb E^n$.
	%
	%
The framework $G(p)$ is 
\emph{first-order rigid} in $\bb S_+^n$ 
if all first-order motions of $G(p)$ are restrictions
of trivial first-order motions.

\subsection{Remark} Note that the equations in (\ref{f.o.m in S}) are equivalent
to the following conditions,
	$$
		(p_i - p_j) \cdot (u_i - u_j) = 0 \and p_k \cdot u_k = 0,
	$$
which are similar to the equations defining first-order rigidity in
$\bb E^n$.

\subsection{Remark} If $G(p)$ is a bar-and-joint framework in $\bb S_+^n$, then 
the graph obtained from $G$ by adjoining a new vertex with edges incident 
with all vertices of $G$, together with the map 
$\hat p: V \cup \{v+1\} \to \bb E^{n+1}$ given by  
	$$ \hat p(i) = \begin{cases}
		p(i) & \text{if } i \neq v + 1 \\
		0 & \text{if } i = v+1 \\
		\end{cases},$$
is first-order rigid in $\bb E^{n+1}$ iff $G(p)$ is first-order rigid in
$\bb S_+^{n+1}$. That is, frameworks in $\bb S_+^n$ can be modeled by the cone on
the same framework in $\bb E^{n+1}$. 


\section{Equivalence of First-Order Rigidity in $\bb S_+^n$ and $\bb E^n$.}
\label{equivalence S -> E}

This section presents two maps, a map carrying a framework $G(p)$ in
$\bb S_+^n$ into a framework $G(q)$ in $\bb E^n$, and a map carrying the
first-order
motions of $G(p)$ into first-order motions of $G(q)$. The latter
map carries trivial 
first-order motions of $\bb S_+^n$ to trivial first-order motions
of $\bb E^n$, yielding the
result $G(p)$ is first-order rigid iff $G(q)$ is first-order rigid.

\subsection{Mapping frameworks and first-order motions} 
\label{mapping motions}
If $G(p)$ is a framework in $\bb S_+^n$,
then $G(\psi \circ p)$ is a framework in $\bb E^n$, where $\psi : \bb
S^n \rightarrow \bb E^n$ is given by $\psi(x) = x/(\e \cdot x)$.
The inverse of $\psi$ is given by
	$\psi^{-1}(x) = x / \sqrt{ x \cdot x }$.
\begin{figure}[htb]
  \begin{center}
    \includegraphics*{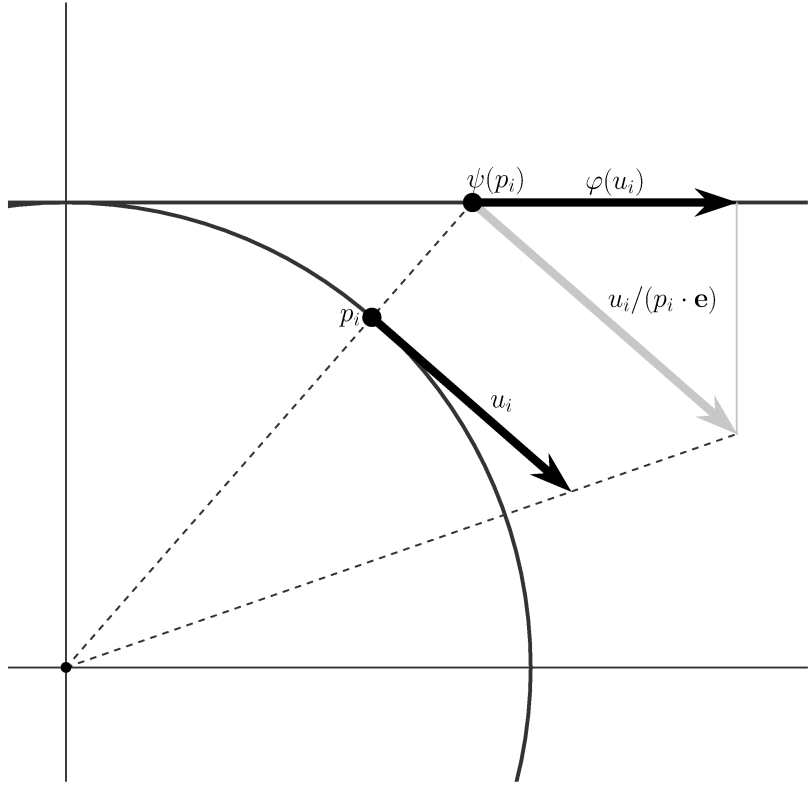}
    \caption{Mapping first-order motions of a framework in $\bb S^n_+$
	to first-order motions of a framework in $\bb E^n$.}
    \label{figure:sphere to plane}
  \end{center}
\end{figure}

If $u$ is a first-order motion of the framework $G(p)$ in $\bb S_+^n$, let
$\varphi$ denote the map
	$$      \varphi : u_i \mapsto 
                \frac{1}{\e \cdot p_i} \left(u_i 
                - (u_i \cdot \e) \e \right).
	$$
If $G(q)$ is a framework
in $\bb E^n$ with first-order motion $v$, then $\varphi^{-1}$ is given
by
	$$\varphi^{-1}: v_i \mapsto
					\frac{1}{\sqrt{q_i \cdot q_i}} \left(
					v_i - (v_i \cdot q_i) \e
					\right).	
	$$
Observe that $\varphi$ and $\varphi^{-1}$ map into the appropriate 
tangent spaces: $\psi^{-1}(q_i) \cdot \varphi^{-1}(v_i) = 0$
and $\varphi(u_i) \cdot \e = 0.$

\subsection{Theorem} \label{theoremSE}
$u$ is a first-order motion of the framework $G(p)$ in $\bb S_+^n$ iff
$\varphi \circ u$ is a first-order motion of the framework $G(\psi \circ p)$
in $\bb E^n$. Moreover, $u$ is a trivial first-order motion iff
$\varphi \circ u \circ \psi^{-1}$ is a trivial first-order motion.

\Pf\ Note that 
	\begin{align} \label{expansion}
	 (\psi(p_i) - \psi(p_j)) \cdot (\varphi(u_i) - \varphi(u_j))
	 = 
	 \frac{p_i \cdot u_i}{(\e \cdot p_i)^2}
	 - \frac{p_i \cdot u_j + p_j \cdot u_i}{(\e \cdot p_i)(\e \cdot p_j)}  
	 +  \frac{p_j \cdot u_j}{(\e \cdot p_j)^2}.
	\end{align}
If $u$ is a first-order motion of $G(p)$, then $u_i \cdot p_i = 0$ 
for all $i \in V$, and $p_i \cdot u_j + p_j \cdot u_i = 0$
for all $\{i,j\} \in E$. By (\ref{expansion}),
$(\psi(p_i) - \psi(p_j)) \cdot (\varphi(u_i) - \varphi(u_j)) = 0$
for all $\{i,j\} \in E$. 
The definition of $\varphi$ ensures that $\varphi(u_i) \cdot \e = 0.$
Therefore, $\varphi \circ u$ is a first-order motion of $G(\psi \circ p)$.

Conversely, 
suppose $\varphi \circ u$ is a first-order motion of $G(\psi \circ p)$. 
Then for all $\{i,j\} \in E$, $(\psi(p_i) - \psi(p_j)) \cdot (\varphi(u_i) - \varphi(u_j)) = 0$.
The observation at the end of the \ref{mapping motions} gives that
$p_i \cdot u_i = \psi^{-1}(\psi(p_i)) \cdot \varphi^{-1}(\varphi(u_i)) = 0$ 
for all $i \in V$. Equation (\ref{expansion}) reduces to 
	$p_i \cdot u_j + p_j \cdot u_i = 0.$
So $u$ is a first-order motion of $G(p)$.

Suppose $u$ is a trivial first-order motion. Then 
$x \cdot u(x) = 0$ for all $x \in \bb S_+^n$ and 
$x \cdot u(y) + y \cdot u(x) = 0$ for all
$x, y \in \bb S_+^n$. Let $v: \bb E^n \to \bb R^{n+1}$ 
denote the composition $\phi \circ u \circ \psi^{-1}$. If $\hat x, \hat y 
\in \bb E^n$ with $x$ denoting $\psi^{-1}(\hat x)$ and $y$ denoting 
$\psi^{-1}(\hat y)$, then (\ref{expansion}) gives
	$$
		(\hat x - \hat y) \cdot (v(\hat x) - v(\hat y)) = 
	 \frac{
		x \cdot u(x)}{(\e \cdot x)^2}
	 - \frac{x \cdot u(y) + y \cdot u(x)}
	 			{(\e \cdot x)(\e \cdot y)}  
	 +  \frac{y \cdot u(y)}{(\e \cdot y)^2} = 0.
	$$
So $v$ is a trivial first-order motion. The converse follows similarly.

\corollary $G(p)$ is first-order rigid in $\bb S_+^n$ iff 
$G(\psi \circ p)$ is first-order rigid in $\bb E^n$.

\subsection{Remark} $\bb S_+^n$ versus $\bb S^n$: Given a discrete
framework, there exists a rotation of the $n$-sphere such that no vertex
of the framework lies on the equator of the sphere. Therefore, we need
not restrict our frameworks to a hemisphere.


\section{Equivalence of First-Order Rigidity in Other Geometries.}
\label{equivalence in other geometries}

\subsection{Geometries.} For 
$x$, $y \in \bb R^{n+1}$, let $\lb x, y \rb_k$ denote
the function
	$$
        \lb x, y \rb_k = x_1y_1 + \cdots + x_{n-k+1}y_{n-k+1} -
        x_{n-k+2}y_{n-k+2} - \cdots - x_{n+1}y_{n+1},
	$$
and let $X_{c,k}^n$ denote the set,
	$$
		X_{c,k}^n = \{ x \in \bb R^{n+1} \mid \lb x, x \rb_k = c, x_{n+1} > 0 \},
	$$ 
for some constant $c \neq 0$ and $k \in \bb N$. We write $X^n$ to simplify
notation, if $c$ and $k$ are understood. If $k=1$ and $c=-1$, then $X^n$ is
\emph{hyperbolic space}, $\bb H^n$.  If $k=1$ and $c=1$, then $X^n$ is
\emph{exterior hyperbolic space}, $\bb D^n$. Spherical space $\bb S_+^n$
is the case $k=0$, $c=1$. Note that $\bb E^n \neq X^n$ for
any choice of $c$ and $k$.

\subsection{Remark} In more generality we can replace $\lb x, y \rb_k$ 
with
	$$
        \lb x, y \rb = a_1x_1y_1 + \cdots + a_{n+1}x_{n+1}y_{n+1},
	$$
 where $a_i \neq 0$ for all $i$, with the exception for Euclidean space:
$a_1 = a_2 = \cdots = a_n = 1$ and $a_{n+1} = 0$.

\subsection{First-order rigidity in $X^n$} 
\label{first-order rigidity in X}
 A metric $d_X$ can be placed on $X^n$ so that $d_X(x,y)$ is a function
of $\lb x, y \rb_k$. A sufficient condition for the distance $d_X(x,y)$
remaining constant is the requirement $\lb x, y \rb_k$ remain constant.
Therefore, the same analysis motivates the following extensions of the
definitions of first-order rigidity to $X^n$.

A \emph{bar-and-joint framework} $G(p)$ in $X^n$ 
is a graph $G$ together with a map $p: V \rightarrow X^n$. 
A \emph{first-order motion} of the
framework $G(p)$ in $X^n$ is a map $u: V \rightarrow \bb R^{n+1}$ satisfying
for each $\{i,j\} \in E$, 
\begin{equation}  \label{eqn_oneX}
        \lb p_i, u_j \rb_k + \lb p_j, u_i \rb_k = 0,
\end{equation}
and for each $i \in V$,
\begin{equation} \label{eqn_twoX}
        \lb p_i, u_i \rb_k = 0.
\end{equation}
A \emph{trivial first-order motion} of $X^n$ is a map $u: X^n \to \bb R^{n+1}$
satisfying
	$$
	        \lb x, u(y) \rb_k + \lb y, u(x) \rb_k = 0
	        \and
           \lb z, u(z) \rb_k = 0
	$$
for all $x,y,z \in X^n$. $G(p)$ is \emph{first-order rigid} in $X^n$ if
all first-order motions of $G(p)$ are the restrictions of trivial 
first-order motions of $X^n$.

\subsection{$X^n$ and $\bb E^n$.} 
In section \ref{equivalence S -> E} we established the 
equivalence between first-order rigidity in $\bb E^n$ and 
first-order rigidity in $\bb S_+^n$. We need only demonstrate the
equivalence holds between the first-order rigidity theories of 
$X^n$ and $\bb S_+^n$.

\subsection{$X^n$ and $\bb S_+^n$} 
Let $\psi_{\bb S_+}: X^n \rightarrow \bb S_+^n$ denote the map $x \mapsto x /
\sqrt{x \cdot x}$, and let $\varphi_{\bb S_+}$ denote the map
	$$
		\varphi_{\bb S_+}: u_i \mapsto \frac{J_k(u_i)}{\sqrt{p_i \cdot p_i}},
	$$ 
where $J_k(x) = (x_1, \cdots, x_{n-k+1}, -x_{n-k+2}, \cdots, -x_{n+1})$.

\begin{figure}[htb]
  \begin{center}
    \includegraphics[scale=0.40]{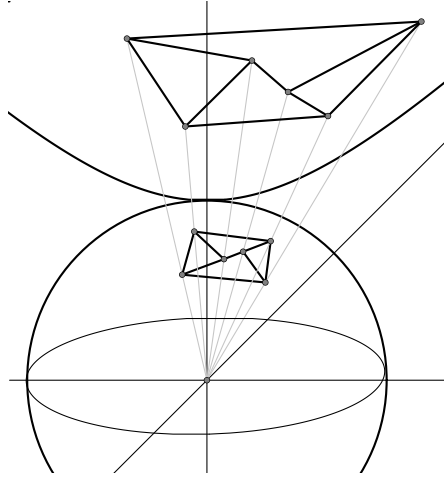}
    \caption{Mapping a bar-and-joint framework from the spherical plane
	 $\bb S_+^2$ into the hyperbolic plane $\bb H^2$.}
    \label{figure: transfering frameworks}
  \end{center}
\end{figure}

\begin{figure}[htb]
  \begin{center}
    \includegraphics{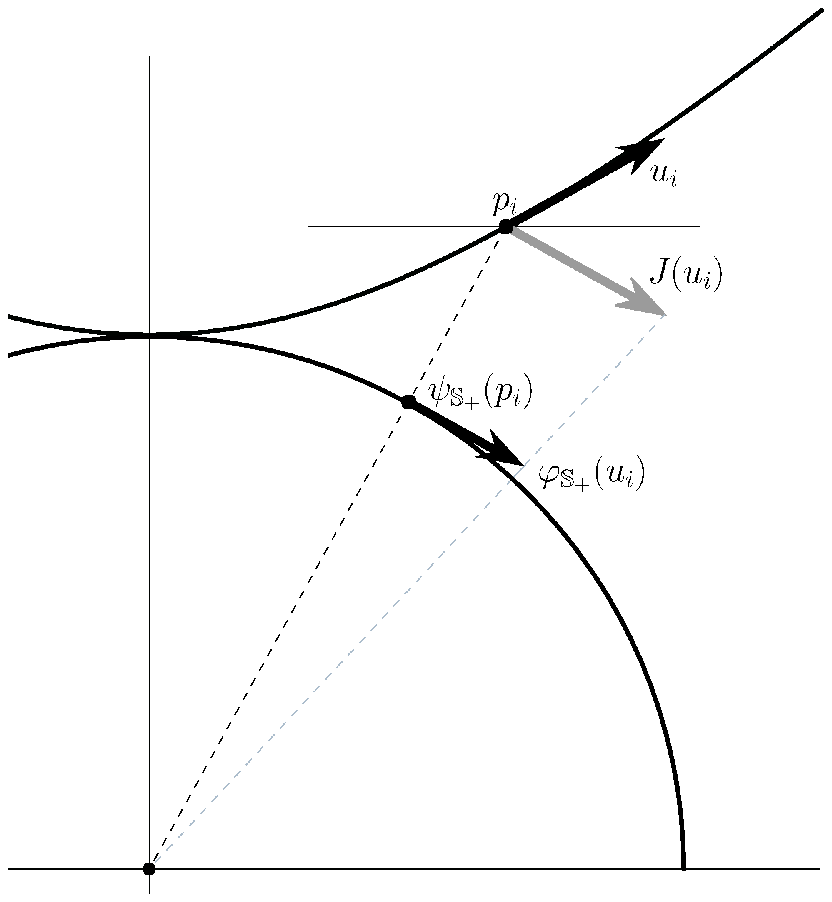}
    \caption{Mapping first-order motions of a framework in $\bb S^n_+$
	to first-order motions of a framework in $\bb H^n$.}
    \label{figure: hyperboloid to sphere}
  \end{center}
\end{figure}

\subsection{Theorem} $G(p)$ is first-order rigid in $X^n$ iff 
$G(\psi_{\bb S_+} \circ p)$ is first-order rigid in $\bb S_+^n$.

\Pf\ Since, $\lb x, y \rb_k = x \cdot J_k(y)$ we have
	\begin{align*}
        \left(\psi_{\bb S_+}(p_i) - \psi_{\bb S_+}(p_j)\right) \cdot 
        \left(\varphi_{\bb S_+}(u_i) - \varphi_{\bb S_+}(u_j)\right) 
         = \frac{\lb p_i, u_i \rb_k}{p_i \cdot p_i} -
        \frac{\lb p_i, u_j \rb_k + \lb p_j, u_i \rb_k}
        {\sqrt{p_i \cdot p_i}\sqrt{p_j \cdot p_j}} +
        \frac{\lb p_j, u_j \rb_k}{p_j \cdot p_j}.
	\end{align*}
 As in the proof of Theorem \ref{theoremSE}, the above equation and the
definitions of $\psi_{\bb S_+}$ and $\varphi_{\bb S_+}$ give that
$\varphi_{\bb S_+} \circ u$ is a first-order motion of $G(\psi_{\bb S_+}
\circ p)$ iff $u$ is a first-order motion of $G(p)$.

It is clear that trivial motions of $\bb S_+^n$ map to trivial motions
of $X^n$. However, a trivial motion of $X^n$ maps onto a ``trivial
motion'' of a proper subset of $\bb S_+^n$. The following fact finishes 
of this proof.

\medskip\noindent\textbf{Fact.} Given a first-order motion $u$ of
$K_{n+1}$, the complete graph on $n+1$ vertices in $\bb E^n$, there
exists a unique trivial first-order motion of $\bb E^n$ extending $u$.

\medskip

(This result and the equivalence of the first-order theories of $\bb
E^n$ and $\bb S_+^n$ give the corresponding result for $\bb S_+^n$,
which was needed to finish the proof of the proceeding theorem.)
\begin{figure}[htb]
  \begin{center}
    \includegraphics[scale=1.25]{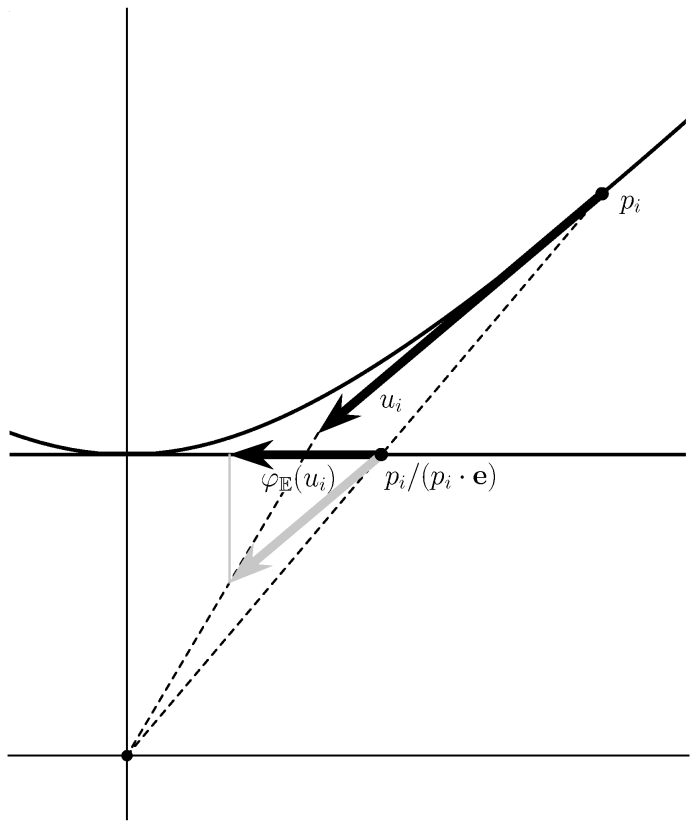}
    \caption{Mapping first-order motions of a framework in $\bb S^n_+$
	to first-order motions of a framework in $\bb E^n$.}
    \label{figure: hyperbolic to plane}
  \end{center}
\end{figure}

\subsection{Remark} There is no obstruction to defining a framework with
vertices in $\bb H^n$ and $\bb D^n$: the equations defining first-order
motions provide formal constraints between these vertices, although the
geometric interpretations of these constraints may not be obvious. In
general, the theorem holds for frameworks with vertices on the surface
$\lb x, x \rb_k = \pm 1$, but not with vertices on $\lb x, x \rb_k = 0$.


\section{The Rigidity Matrix}
\label{rigidity matrix}

\subsection{Projective models of $X^n$} The projective model of $X^n$ is 
the subset of $\bb E^n$ obtained by projecting from the origin 
the points of $X^n$ onto $\bb E^n$,
	$$
		\left\{ \frac{1}{\e \cdot x}\ x \ \Big| \
							x \in X^n \right\} \subset \bb E^n.
	$$
 The projective model of hyperbolic $n$-space $\bb H^n$ is the interior
of the unit $n$-ball $B^n$ of $\bb E^n$ and the projective model of
exterior hyperbolic $n$-space $\bb D^n$ is the exterior of $B^n$. The
unit $(n-1)$-sphere $S^{n-1}$ is the \emph{absolute}, the points at
infinity of hyperbolic geometry. Spherical $n$-space is model
projectively by $\bb E^n$. 

Since we are now restricting our attention to points in $\bb E^n$, we
identify $\bb E^n$ with $\bb R^n$ and write $PX^n$ to denote the
projective model of $X^n$ as a subset of $\bb R^n$. Distance in $PX^n$ 
is calculated by normalizing the points
into $X^n$ and applying the definition of distance in $X^n$. For example,
the distance between points $x$ and $y$ in $P\bb S_+^n$ (so $x,y \in \bb R^n$)
is 
	$$
	d_{P\bb S_+}(x,y) = \arccos \left(\frac{ 1 + x \cdot y }
        {\sqrt{ 1 + x \cdot x} \sqrt{ 1 + y \cdot y}}\right),
	$$
and for points $x$ and $y$ in $P\bb H^n$,
	$$
	d_{P\bb H}(x,y) = \text{arccosh} \left(\frac{ 1 - x \cdot y }
        {\sqrt{ 1 - x \cdot x} \sqrt{ 1 - y \cdot y}}\right).
	$$

\subsection{The rigidity matrix of a framework}

A first-order motion $u: V \to \bb R^n$ of the framework $G(p)$ in $\bb R^n$,
satisfies 
	\begin{equation*}\label{equation euclidean fom}
		(p_i - p_j) \cdot (u_i - u_j) = 0.
	\end{equation*}
This system of homogeneous linear equations, indexed by the edges of $G$,
induces a linear transformation with matrix $R_{\bb E}(G,p)$, called
the \emph{rigidity matrix} of $G(p)$,
	$$
	R_{\bb E}(G,p) = \bordermatrix{ 
	        &        & i           & \cdots & j           &        & \cr
	        &        & \vdots      &        & \vdots      &        & \cr
	\{i,j\} & \cdots & p_i - p_j & \cdots & p_j - p_i & \cdots & \cr
	        &        & \vdots      &        &  \vdots     &        & \cr}.
	$$
The kernel of $R_{\bb E}(G,p)$ is precisely the space of first-order motions
of $G(p)$.

A first-order motion $u: V \to \bb R^n$ of the framework $G(p)$ in $P\bb H^n$,
$P\bb D^n$ or $P\bb S^n$ satisfies
	\begin{equation*}\label{equation non-euclidean fom}
	(k_{ij} + k_{ji}) \cdot (u_i + u_j) = 0,
	\end{equation*}
where $k_{ij}$ is
	$$
	k_{ij} =	
	\begin{cases}
	\left(\frac{1 - p_i \cdot p_j}{1 - p_i \cdot p_i}\right) p_i - p_j, 
		& \text{ for } P\bb H^n \text{ or } P\bb D^n \\
	\left(\frac{1 + p_i \cdot p_j}{1 + p_i \cdot p_i}\right) p_i - p_j, 
		& \text{ for } P\bb S_+^n
	\end{cases}.
	$$
The matrix of the linear transformation induced by this system of linear
equations is the \emph{rigidity matrix} $R_{X}(G,p)$ of $G(p)$,
	$$
	R_{X}(G,p) = \bordermatrix{ 
	        &        & i           & \cdots & j           &        & \cr
	        &        & \vdots      &        & \vdots      &        & \cr
	\{i,j\} & \cdots & k_{ij} & \cdots & k_{ji} & \cdots & \cr
	        &        & \vdots      &        &  \vdots     &        & \cr}.
	$$
Note that $k_{ij}$ depends on $X$.

\begin{figure}[htb]
  \begin{center}
     \includegraphics{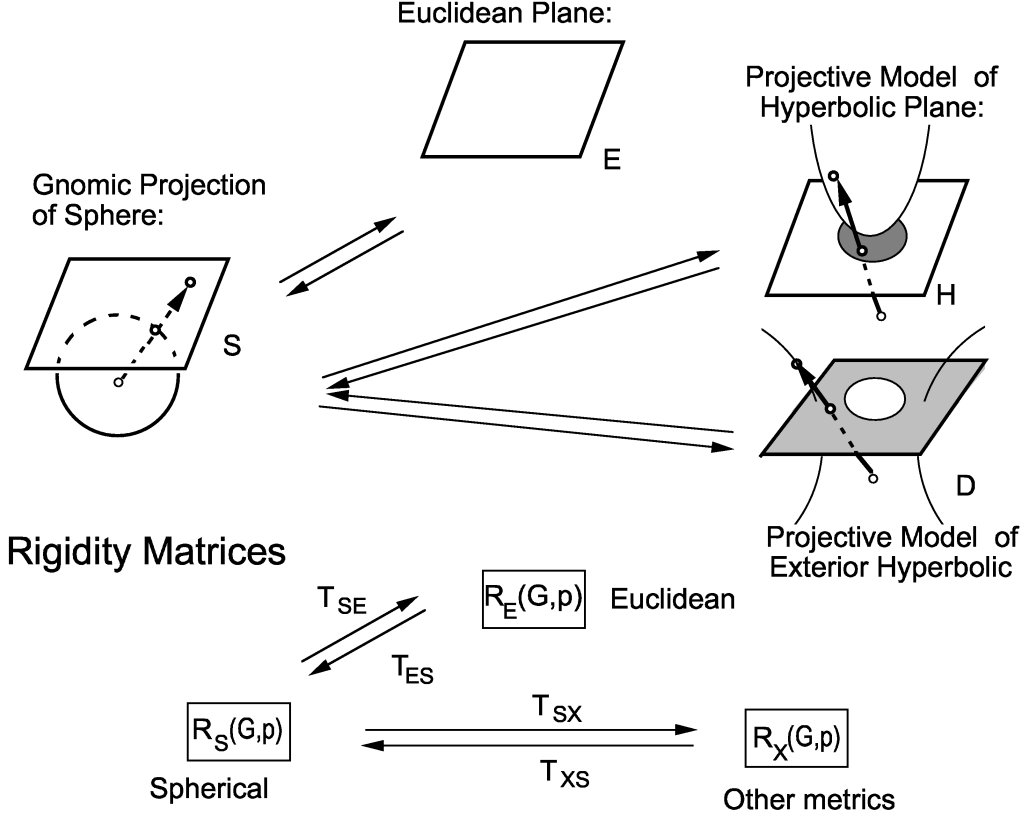}
    \caption{A visual summary of the equivalence of first-order rigidity
		in the projective models of hyperbolic geometry $H$,
		spherical geometry $S$ and Euclidean geometry $E$. Here $T_{SE}$
		denotes the linear transformation $T_1(G,p)$ defined in the text, 
		$T_{ES}$ the inverse of $T_{SE}$.}
    \label{figure: summary of theorem}
  \end{center}
\end{figure}

\subsection{Transforming rigidity matrices} Let $T_K(G,p)$ denote the matrix
	$$
	T_K(G,p) = 
		\left[\begin{array}{cccc}
		T_{p_1} & 0 & 0 & 0 \\
		0 & T_{p_2} & 0 & 0 \\
		0 & 0 & \ddots & 0 \\
		0 & 0 & 0 & T_{p_{v}} 
		\end{array}\right],
	$$
where $T_{p_k} = I + K (p_k^{(i)}p_k^{(j)})$ ($I$ is the $n
\times n$ identity matrix and $(p_k^{(i)}p_k^{(j)})$ is the $n \times n$
matrix with $p_k^{(i)}p_k^{(j)}$ as entry $(i,j)$, where $p_k^{(i)}$ is
the $i$-th component of $p_k$). For example, for 
$n=3$ and $p_{k} = (x_1, x_2, x_3)$,
	$$
	T_{p_k} = 
		\left[\begin{array}{ccc}
		1 + Kx_1^2 & Kx_1x_2 & Kx_1x_3 \\
		Kx_1x_2 & 1 + Kx_2^2 & Kx_2x_3 \\
		Kx_1x_3 & Kx_2x_3 & 1 + Kx_3^2
		\end{array}\right].
	$$

\theorem Let $G(p)$ be a framework with $p \in \bb R^n$. 
Then 
\begin{enumerate}
\item
$T_K(G,p)$ satisfies
	\begin{align*}
	R_{P\bb H} \times T_{-1}(G,p) = R_{\bb E}(G,p)
		\and 
	R_{P\bb S_+} \times T_{1}(G,p) = R_{\bb E}(G,p);
	\end{align*} 
\item 
$G(p)$ is first-order rigid in $P\bb S_+^n$ iff
$G(p)$ is first-order rigid in $P\bb E^n;$
\item
$G(p)$ is first-order rigid in $P\bb H^n \cup P\bb D^n$ iff
$G(p)$ is first-order rigid in $P\bb E^n$ and $p_i \cdot p_i \neq 1$
for all $i \in V$ (no vertex is on the absolute).
\end{enumerate}

\Pf\ (1) Since $T_{p_i}$ multiplies only the columns corresponding to vertex
$i$, we need only verify $k_{ij} \times T_{p_i} = p_i - p_j$. 
This is a straightforward calculation, 
\begin{eqnarray*}
& & 
k_{ij} \times \left(\text{column } \ell \text{ of } T_{p_i} \right) \\
& = &  
	\left(\frac{1 + K (p_i \cdot p_j) }{1 + K (p_i \cdot p_i)} 
	p_i - p_j\right) \cdot \left(\e_{\ell} + K p_i^{(\ell)}
	 p_i\right)
\\
& = &  
	\left(\frac{1 + K (p_i \cdot p_j) }{1 + K (p_i \cdot p_i)}\right)
	\left( p_i \cdot \e_{\ell} + K p_i^{(\ell)} (p_i \cdot p_i)\right)
	- \left( p_j^{(\ell)} + K p_i^{(\ell)} (p_j \cdot p_i)	\right)
\\
& = & 
	\left(\frac{1 + K (p_i \cdot p_j) }{1 + K (p_i \cdot p_i)}\right)
	\left( 1 + K (p_i \cdot p_i) \right) p_i^{(\ell)}
	- 	\left( p_j^{(\ell)} + K p_i^{(\ell)} (p_j \cdot p_i)	\right)
\\
& = & 
	\left( {1 + K (p_i \cdot p_j) } \right) p_i^{(\ell)}
	- 	\left( p_j^{(\ell)} + K p_i^{(\ell)} (p_j \cdot p_i)	\right)
\\
& = & 
	p_i^{(\ell)} + K p_i^{(\ell)} (p_i \cdot p_j) 
	- p_j^{(\ell)} - K p_i^{(\ell)} (p_j \cdot p_i)
\\
& = & 
	p_i^{(\ell)} - p_j^{(\ell)},
\end{eqnarray*}
which is column $\ell$ of $p_i - p_j$.

(2), (3): Since the determinant of $T_K(G,p)$ is the product 
$\prod_{i=1}^v \det(T_{p_i})$ and 
 	$$
	\det(T_{p_i}) = 1 + K(p_i \cdot p_i),
	$$
 the dimension of the vector space of first-order motions of $G(p)$ 
is the same in each geometry iff $1 + K(p_i \cdot p_i) \neq 0$ for 
all $i \in V$.

\subsection{Remark} It is well-known that the rank of the rigidity
matrix, and thus first-order rigidity, of a framework in $\bb E^n$ is
invariant under projective transformations of $\bb E^n$. Due to the
equivalence of first-order theories, the same is true of frameworks in
$X^n$. (In fact, there exists an underlying projective theory.)
Intuitively at least, this projective invariance suggests the
equivalences presented in this paper since all the geometries discusses
can be obtained from projective geometry by choosing an appropriate set
of transformations.
\begin{figure}[htb]
  \begin{center}
    \includegraphics[trim=100 275 100 275]{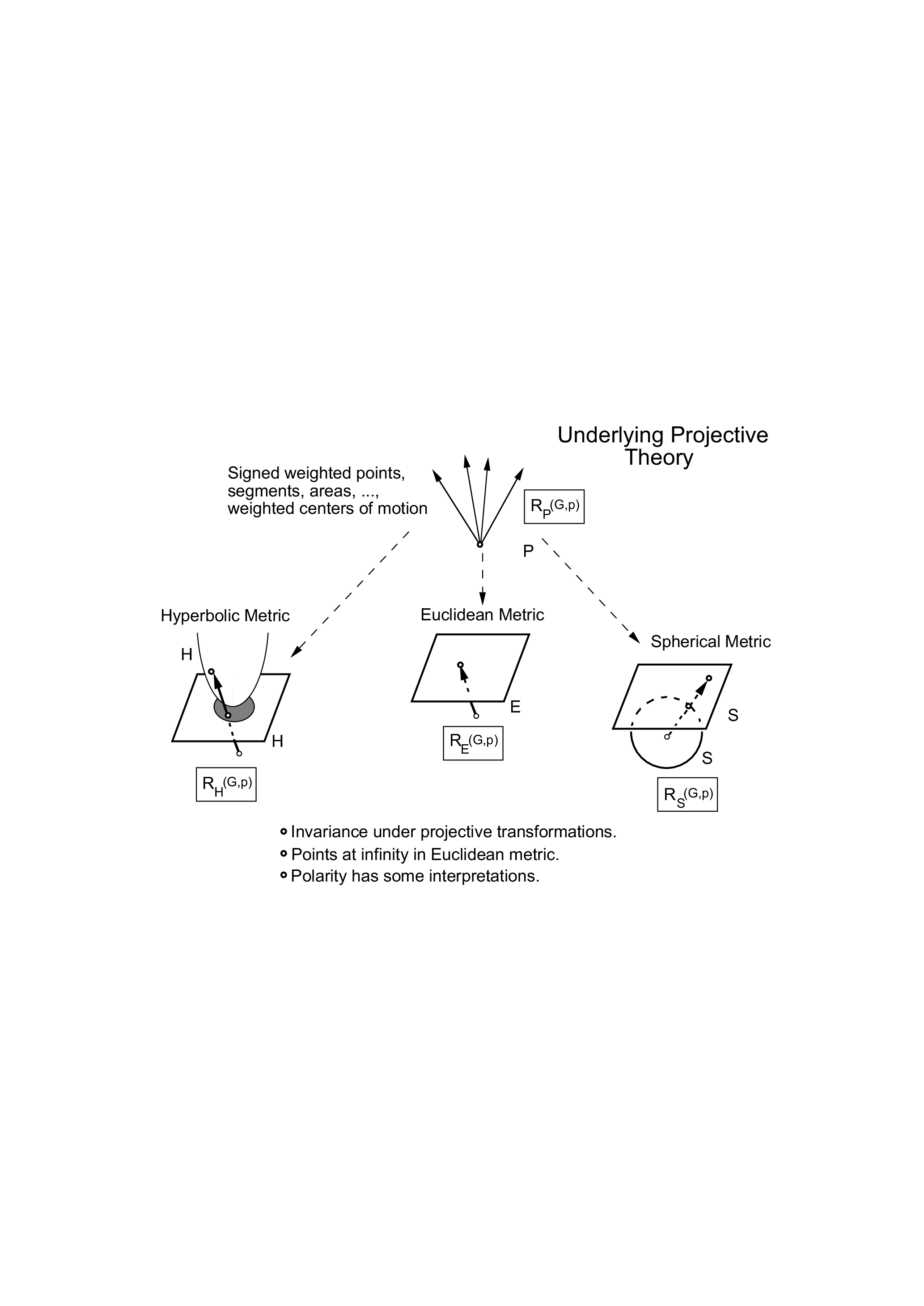}
    \caption{A visual summary of the 
		underlying projective theory: hyperbolic space $H$, Euclidean space $E$
		and spherical space $S$ can be realized as \emph{subgeometries} of
		projective geometry.}
    \label{figure: summary of projective theory}
  \end{center}
\end{figure}


\section{The First-Order Uniqueness Theorems of Andreev and Cauchy-Dehn}
\label{andreev}

An immediate consequence of the equivalence of these first-order rigidity
theories is the ability to transfer results between the theories.



\subsection{The Cauchy-Dehn Theorem} 

The Cauchy-Dehn theorem for polytopes in $\bb E^n$ states that a convex,
triangulated polyhedron in $\bb E^n$, $n \geq 3$, is first-order rigid.
Before the generalization of this theorem can be stated, convexity in $X^n$
needs to be defined. A set $S \subset X^n$ is \emph{convex} if, for any
line $L$ of $X^n$, $L \cap S$ is connected. Therefore, $S \subset X^n$
is convex iff $\psi_{\bb E}(S) \subset \bb E^n$ is convex. 

\theorem (Cauchy-Dehn)
  A convex, triangulated polytope $P$ in $X^n$, $n \geq 3$,
  is first-order rigid.

\subsection{A first-order version of Andreev's uniqueness theorem}
If $p$ denotes a point of $\bb D^n$, then the set of points $x$ in 
$\bb R^{n+1}$ satisfying $\lb p, x \rb_{1} = 0$ (orthogonal in the
hyperbolic sense) defines a unique hyperplane of $\bb R^{n+1}$ through
the origin. Therefore, to each point of $p$, there corresponds a 
unique hyperplane of $\bb H^n$,
	$$
	P = \{ x \in \bb H^n \mid \lb p, x \rb_1 = 0 \},
	$$ 
and conversely.

If $q$ is another point of $\bb D^n$ with $Q$ the corresponding hyperplane
of $\bb H^n$, the angle of intersection of the hyperplanes $P$ and $Q$
is defined to be $\arccos( \lb p, q \rb_1 )$. So 
equations~(\ref{eqn_oneX}) and ~(\ref{eqn_twoX}) 
defining a first-order motion $u$ of a framework $G(p)$ 
in $\bb D^n$, 
	$$
        \lb p_i, u_j \rb_k + \lb p_j, u_i \rb_k = 0
			\and
        \lb p_i, u_i \rb_k = 0,
	$$
are precisely the conditions defining a ``first-order motion''
of a collection of planes under angle constraints (a bar-and-joint 
framework is merely a collection of points under distance constraints). 
Polyhedra with fixed dihedral angles are examples of such objects.

Under this point-plane correspondence of $\bb D^n$ and $\bb H^n$, the
Cauchy-Dehn theorem for $\bb D^n$ gives a first-order version of Andreev's
uniqueness theorem. Indeed, a simple, convex polytope in $\bb H^n$ is a
triangulated, convex polytope in $\bb D^n$. We use \emph{stiff} to
denote the analogous definition of first-order rigid.

\theorem (Andreev)
  If $M$ is a simple, convex polytope in $\bb H^n$, $n \geq 3$,
  then $M$ is stiff.

\subsection{Remark} The usual hypothesis of Andreev's theorem requires the
polytope $M$ to have dihedral angles not exceeding $\pi/2$. This 
supposition implies $M$ is simple.

\subsection{Remark} The point-plane correspondence described above is
known as polarity. There is a version of this result for the sphere that
requires a better discussion of polarity on the sphere.




  \end{document}